\documentclass[10pt,twoside,reqno]{amsart}

\usepackage[T1]{fontenc}
\usepackage[utf8]{inputenc}
\usepackage[english]{babel}
\usepackage{lmodern}

\usepackage{amsmath,amsfonts,amsthm,amssymb,amscd}
\usepackage{mathrsfs}
\usepackage{latexsym}
\usepackage{graphicx}
\usepackage{enumitem}

\usepackage{geometry}
\theoremstyle{plain}
\newtheorem{theorem}{Theorem}
\newtheorem{lemma}{Lemma}
\newtheorem{corollary}{Corollary}

\theoremstyle{remark}

\usepackage{hyperref}
\hypersetup{
  colorlinks=true,
  linktoc=all,
  linkcolor=blue,
  citecolor=blue,
  urlcolor=blue,
  bookmarksopen=true
}
\AtBeginDocument{%
 }

\allowdisplaybreaks[2]

\title[Mean estimates for short sums]{Mean Estimates for Short Polynomial Exponential Sums over Primes}

\author{Karimjon Ibrohimjonovich Mirzoabdughafurov}
\date{}

\address{\normalfont\bf Tajik State University of Finance and Economics, Dushanbe, Tajikistan}

\email{\bf karimjon2003@mail.ru}
\subjclass[2020]{Primary 11L20; Secondary 11L07, 11L15}
\keywords{short exponential sums, prime numbers, mean estimates, Fej\'er kernel, Weyl's method, rational approximations}

\begin{document}
\begin{abstract}
Let $n\ge2$ be a fixed integer, and let $K,x,y$ be positive
integers satisfying $2\le K\le y<x$. Consider a polynomial
\[
f(u)=\alpha u^n+\alpha_{n-1}u^{n-1}+\cdots+\alpha_1u+\alpha_0
\]
with real coefficients. We study the quantity
\[
V_n(K,x,y)
=
\sum_{k=1}^{K}
\left|
\sum_{x-y<p\le x}e(kf(p))
\right|,
\qquad e(t)=e^{2\pi it},
\]
where the inner sum is over primes. For 
$\alpha=\frac{a}{q}+\frac{\theta}{q^2}$,
$a\in\mathbb Z$, $q\in\mathbb N$,
$(a,q)=1$, $|\theta|\le1$,
we establish the bound 
\[
\begin{aligned}
V_n(K,x,y)
&\ll Ky\Biggl[
\frac{1}{\sqrt{K\log(2y)}}+
\min\left\{
\Delta^{\frac{1}{2^n}}
\bigl(\log(2y)\bigr)^{\frac{n^2-1}{2^n}},
\Delta^{\frac{1}{3\cdot 2^{n-2}}}
\bigl(\log(2y)\bigr)^{\frac{n^3-1}{3\cdot 2^{n-1}}}
\right\}
\Biggr],
\end{aligned}
\]
where 
\[
\Delta=
\frac{1}{q}+\frac{1}{y}+\frac{q}{Ky^n}.
\]
The bound is uniform in the position of the interval, and the logarithmic factors depend on its length. The proof uses the nonnegative Fej\'er kernel, successive differencing, and the second and third moments of the generalized divisor function. The Brun--Titchmarsh inequality accounts for the number of primes in the interval and provides an additional logarithmic saving in the term arising from averaging. We obtain sufficient conditions for a saving of any fixed power of the logarithm of the interval length, as well as an analogous mean estimate for sums over all integers.
\end{abstract}
\maketitle

Estimates for exponential sums with polynomial phases play a central role in analytic number theory. The classical differencing method goes back to H.~Weyl~\cite{Weyl1916}. Successive differencing reduces the degree of the polynomial phase and hence reduces the estimation of the original sum to estimates for linear sums.

I.~M.~Vinogradov's method of exponential sums is one of the principal tools for studying additive problems and the distribution of arithmetic sequences. Its development is presented in Vinogradov's monograph~\cite{Vinogradov1980} and the survey by I.~M.~Vinogradov and A.~A.~Karatsuba~\cite{VinogradovKaratsuba1984}. Estimates for Weyl sums and differencing methods are also discussed in the monographs of N.~M.~Korobov~\cite{Korobov1989}, R.~Vaughan~\cite{Von}, and H.~Montgomery~\cite[Chapters~3--4]{Montgomery1994}.

Short exponential sums over primes arise when the Hardy--Littlewood method is applied to additive problems with restrictions on the sizes of the summands. In particular, they occur in the Waring--Goldbach problem with almost equal prime summands. In~\cite{WeiWooley2015}, Weyl-type estimates for sums over short intervals are applied to representations of integers as sums of powers of almost equal primes.

Pointwise estimates for short Weyl sums in an intermediate range of rational approximation and their applications to Waring's problem with almost equal and almost proportional summands were studied by the author in~\cite{Mirzoabdughafurov2026Quartic,Mirzoabdughafurov2026PPS}.

Alongside estimates for individual exponential sums, it is of interest to estimate the sum of their absolute values over an integer multiplier of the phase. For short intervals, the explicit dependence of such estimates on the interval length, the averaging parameter, and the quality of the rational approximation to the phase coefficient is essential.

Let $n\ge2$ be a fixed integer, and let $K$, $x$, $y$ be positive
integers satisfying $2\le K\le y<x$. Consider a polynomial
\[
f(u)=\alpha u^n+\alpha_{n-1}u^{n-1}+\cdots+\alpha_1u+\alpha_0
\]
with real coefficients, whose leading coefficient satisfies
\begin{equation*}
\alpha=\frac{a}{q}+\frac{\theta}{q^2},
\qquad
a\in\mathbb Z,\quad q\in\mathbb N,\quad
(a,q)=1,\quad |\theta|\le1.
\end{equation*}
For this fixed polynomial, put
\[
V_n(K,x,y)
=
\sum_{k=1}^{K}
\left|
\sum_{x-y<p\le x}e(kf(p))
\right|,
\qquad
e(t)=e^{2\pi it},
\]
where the inner sum is over primes $p$.

For brevity, put 
\[
\Delta=
\frac{1}{q}+\frac{1}{y}+\frac{q}{Ky^n}.
\]

Modern developments in the theory of mean values of Weyl sums are closely connected with the proof of the main conjecture in Vinogradov's mean value theorem. In~\cite{BourgainDemeterGuth2016}, this conjecture was proved for degrees greater than three by the decoupling method. In~\cite{Wooley2019}, the method of nested efficient congruencing was developed, proving the main conjecture for all degrees and allowing generalizations to other polynomial systems. These results underpin modern estimates for Weyl sums and their discrete moments.

Mean estimates for Weyl sums over an integer multiplier were studied by K.~Halupczok~\cite{Halupczok2020}. After shifting the interval and applying the result to the phase under consideration, Theorem~9 of that paper yields 
\[
\sum_{k=1}^{K}
\left|
\sum_{x-y<m\le x}e(kf(m))
\right|
\ll
Ky(Ky)^\varepsilon
\Delta^{\frac{1}{n(n-1)}},
\]
where the bound holds for every fixed 
$\varepsilon>0$, with a constant depending only on $n$ and $\varepsilon$.

Thus, for sums over all integers, bounds are known that give a strong power saving in the Diophantine parameter together with an additional factor of the form $(Ky)^\varepsilon$. In the present paper, a mean estimate over primes is combined with an estimate for sums over integers involving an explicit logarithmic factor.

Sums of short exponential sums over primes are studied by Z.~Kh.~Rakhmonov and F.~Z.~Rakhmonov in~\cite{RakhmonovRakhmonov2014}. They consider estimates for individual short cubic sums weighted by the von Mangoldt function in~\cite{RakhmonovRakhmonov2017}. Mean estimates for short exponential sums over primes are also investigated by F.~Z.~Rakhmonov in~\cite{Rakhmonov2025,Rakhmonov2026}. For the monomial phase $f(u)=\alpha u^n$,
their results yield the bound 
\[
V_n(K,x,y)
\ll 
Ky
\left(
\Delta+K^{-2^{n-1}}
\right)^{\frac{1}{2^{n+1}}}
(\log x)^{\frac{n^2}{2^{n+1}}}.
\]
After separating the terms, this can be written as 
\[
V_n(K,x,y)
\ll 
Ky\left(
\Delta^{\frac{1}{2^{n+1}}}
+
K^{-\frac{1}{4}}
\right)
(\log x)^{\frac{n^2}{2^{n+1}}}.
\]
In the present paper, we obtain two forms of a bound for $V_n(K,x,y)$. In the first, the exponent of $\Delta$
 increases from $\frac{1}{2^{n+1}}$
 to $\frac{1}{2^n}$, and in the second to $\frac{1}{3\cdot 2^{n-2}}$. The corresponding logarithmic factors are 
\[
\bigl(\log(2y)\bigr)^{\frac{n^2-1}{2^n}}
\qquad\text{and}\qquad
\bigl(\log(2y)\bigr)^{\frac{n^3-1}{3\cdot 2^{n-1}}}.
\]
The main theorem combines these two bounds by taking the minimum of the two expressions.

The choice between them depends on the size of $\Delta$
 and the interval length. The second expression does not exceed the first if and only if 
\[
\Delta\bigl(\log(2y)\bigr)^{(n-1)^2(2n+1)}
\le 1.
\]
Thus, the larger exponent of $\Delta$ in the second bound is accompanied by a larger power of the logarithm.

At the same time, the term 
\[
K^{-\frac{1}{4}}
(\log x)^{\frac{n^2}{2^{n+1}}}
\]
is replaced by 
\[
\frac{1}{\sqrt{K\log(2y)}}.
\]
All logarithmic factors in the new bound are expressed in terms of the interval length $y$. The bound is uniform in the position of the interval for $2\le K\le y<x$.

Comparison with the previous bound should take account of both the power and logarithmic factors. In particular, under the condition 
\[
\left(
\Delta+K^{-2^{n-1}}
\right)(\log x)^{n^2}=o(1),
\]
which means that the previous right-hand side is $o(Ky)$, the right-hand side of the bound in Theorem~\ref{theo1} is of smaller order than the previous one.

The main step of the proof is to apply the nonnegative Fej\'er kernel to the mean square of the exponential sums. After expanding the squares, the positivity of the kernel allows summation over one of the two prime variables to be extended to all integers in the interval. Estimating the resulting mixed moment reduces to estimating moments of sums over integers. Retaining the number of primes at this stage allows us to apply the Brun--Titchmarsh inequality~\cite[p.~121, Theorem~2]{MontgomeryVaughan1973} and obtain an additional factor 
$(\log(2y))^{-\frac{1}{2}}$
 in the term arising from averaging.

The moments of sums over integers are estimated by successive differencing. Grouping products of the differencing parameters introduces the divisor function. Applying the Cauchy--Schwarz inequality to the resulting sum and using the second moment of this function gives the first bound. Applying H\"older's inequality and using its third moment gives the second. In both cases, the arithmetic factors are accounted for by fixed powers of the logarithm.

The main theorem yields explicit sufficient conditions on $K$, $y$ and $q$ under which, for a fixed $A>0$, we have 
\[
V_n(K,x,y)
\ll
\frac{Ky}{(\log(2y))^A}.
\]
For comparison, direct estimation of the absolute values and the Brun--Titchmarsh inequality give 
\[
V_n(K,x,y)
\le
K\bigl(\pi(x)-\pi(x-y)\bigr)
\ll
\frac{Ky}{\log(2y)}.
\]
Thus, for $A>1$, the first corollary improves this general upper bound.

The moment estimates also yield a corollary for the inner sum over all integers in the interval. In this case, the term associated with the number of primes is absent, and the result contains the minimum of the same two expressions with explicit logarithmic factors.

For $n=3$ and $n=4$, the exponent 
$\frac{1}{3\cdot 2^{n-2}}$
 in the second expression coincides with the exponent $\frac{1}{n(n-1)}$
 in the bound of K.~Halupczok~\cite{Halupczok2020} quoted above. The additional factor in our result is a fixed power of $\log(2y)$. For other degrees, comparing the results requires taking account of both the exponent of $\Delta$
 and the additional factors.

The following theorem holds.

\begin{theorem}\label{theo1}
Let $n\ge2$ be a fixed integer, let $K,x,y$ be integers
satisfying $2\le K\le y<x$, and let
\[
f(u)=\alpha u^n+\alpha_{n-1}u^{n-1}+\cdots+\alpha_1u+\alpha_0
\]
be a polynomial with real coefficients whose leading coefficient
satisfies
\[
\alpha=\frac{a}{q}+\frac{\theta}{q^2},
\qquad
a\in\mathbb Z,\quad q\in\mathbb N,\quad
(a,q)=1,\quad |\theta|\le 1.
\]

Then we have 
\begin{align*}
V_n(K,x,y)
&\ll
Ky\left[
\frac{1}{\sqrt{K\log(2y)}}
+
\min\left\{
\Delta^{\frac{1}{2^n}}
\bigl(\log(2y)\bigr)^{\frac{n^2-1}{2^n}},\Delta^{\frac{1}{3\cdot 2^{n-2}}}
\bigl(\log(2y)\bigr)^{\frac{n^3-1}{3\cdot 2^{n-1}}}
\right\}
\right].
\end{align*}
The implied constant in Vinogradov's notation depends only on $n$.
\end{theorem}
\begin{corollary}\label{cor1}
Suppose that the hypotheses of Theorem~\ref{theo1} hold, and let $A>0$ be a fixed constant. Put 
\[
B=\min\left\{
2^nA+n^2-1,\;
3\cdot 2^{n-2}A+\frac{n^3-1}{2}
\right\}.
\]
If 
\[
K\gg (\log(2y))^{2A-1},
\qquad
y\gg (\log(2y))^B
\]
and 
\[
(\log(2y))^B
\ll q
\ll 
\frac{Ky^n}{(\log(2y))^B},
\]
then 
\[
V_n(K,x,y)
\ll 
\frac{Ky}{(\log(2y))^A},
\]
where the implied constants in Vinogradov's notation may depend only on $n$ and $A$.
\end{corollary}

The following corollary follows from the moment estimate for the sum $T(h)$ established in the proof of Theorem~\ref{theo1}.

\begin{corollary}\label{cor2}
Under the hypotheses of Theorem~\ref{theo1}, we have 
\begin{align*}
\sum_{k=1}^{K}
\left|
\sum_{x-y<m\le x}e(kf(m))
\right|
\ll 
Ky\min\left\{
\Delta^{\frac{1}{2^n}}
\bigl(\log(2y)\bigr)^{\frac{n^2-1}{2^n}},
\Delta^{\frac{1}{3\cdot 2^{n-2}}}
\bigl(\log(2y)\bigr)^{\frac{n^3-1}{3\cdot 2^{n-1}}}
\right\}.
\end{align*}
The implied constant in Vinogradov's notation may depend only on $n$.
\end{corollary}

\section{Auxiliary results}
\begin{lemma}\label{linear-sum}
Let $\alpha$ be any real number, let 
$Q$ be an integer, and let $P$ be a positive integer. Then 
\[
\left|
\sum_{x=Q+1}^{Q+P} e^{2\pi i\alpha x}
\right|
\le
\min\left(P,\frac{1}{2\|\alpha\|}\right),
\]
where $\|\alpha\|$ is the distance from $\alpha$ to the nearest integer.
\end{lemma}
{\large\bf Proof.} See~{\rm \cite[Lemma~1, p.~12]{Korobov1989}}.
\begin{lemma}\label{Weyl}
For a function $\varphi$, put 
\[
\Delta_1(\varphi(x);h)=\varphi(x+h)-\varphi(x)
\]
and for $j\ge1$ define:
\[
\Delta_{j+1}(\varphi(x);h_1,\ldots,h_{j+1})=\Delta_j(\varphi(x+h_{j+1});h_1,\ldots,h_j)-\Delta_j(\varphi(x);h_1,\ldots,h_j).
\]
Let 
\[
T(\varphi)=\sum_{x=1}^{Q}e\bigl(\varphi(x)\bigr),
\]
where $\varphi$ is any arithmetic function. Then 
\[
|T(\varphi)|^{2^j}
\le
(2Q)^{2^j-j-1}
\sum_{|h_1|<Q}\cdots\sum_{|h_j|<Q}T_j,
\]
where 
\[
T_j=
\sum_{x\in I_j}
e\bigl(\Delta_j(\varphi(x);h_1,\ldots,h_j)\bigr),
\]
and the intervals $I_j=I_j(h_1,\ldots,h_j)$, which may be empty, satisfy 
\[
I_1(h_1)\subset[1,Q],
\qquad
I_j(h_1,\ldots,h_j)
\subset
I_{j-1}(h_1,\ldots,h_{j-1}).
\]
\end{lemma}
{\large\bf Proof.} See~{\rm \cite[Lemma~2.3, p.~18]{Von}}.

\section{Proof of Theorem \ref{theo1}}

Let 
\[
 \mathbb{I}=(x-y,x]\cap\mathbb Z,
\qquad
 \mathbb{P}=\{p\in \mathbb{I}:p\text{ prime}\}.
\]
For $h\in\mathbb Z$, introduce the sums 
\[
S(h)=\sum_{p\in \mathbb{P}}e(hf(p)),
\qquad
T(h)=\sum_{m\in \mathbb{I}}e(hf(m)).
\]
By the Cauchy--Schwarz inequality, 
\begin{equation}\label{cauchy}
V_n^2(K,x,y)
=
\left(\sum_{k=1}^{K}|S(k)|\right)^2
\le K\sum_{k=1}^{K}|S(k)|^2.
\end{equation}

For an integer $M\ge1$, consider the Fej\'er kernel 
\begin{align}
  &\Phi_M(t)
=
\frac{1}{M}\left|\sum_{j=0}^{M-1}e(jt)\right|^2=\frac{1}{M}
\left(\sum_{j=0}^{M-1}e(jt)\right)
\left(\sum_{k=0}^{M-1}e(-k t)\right)\nonumber\\
&=\frac{1}{M}
\sum_{j=0}^{M-1}\sum_{k=0}^{M-1}
e\bigl((j-k)t\bigr)
 =\frac{1}{M}\sum_{j=0}^{M-1}
\sum_{h=j-(M-1)}^{j} e(ht)\nonumber\\
&=
\frac{1}{M}\sum_{h=-(M-1)}^{M-1}
e(ht)
\sum_{j=\max(0,h)}^{\min(M-1,M-1+h)} 1=
\sum_{|h|<M}
\left(1-\frac{|h|}{M}\right)e(ht).\label{Feyer}
\end{align}
To estimate the sum of squares on the right-hand side of~(\ref{cauchy}), apply identity~(\ref{Feyer}) with 
$$
M=2K,\quad \text{and}\quad t=f(p_1)-f(p_2).
$$
Put 
\[
P=\#\mathbb{P}=\pi(x)-\pi(x-y).
\]
Since $S(0)=P$, summing over $p_1,p_2\in \mathbb{P}$, interchanging the finite sums, and taking account of $S(-h)=\overline{S(h)}$, we obtain 
\begin{align}
&\sum_{p_1,p_2\in \mathbb{P}}
\Phi_{2K}\!\left(f(p_1)-f(p_2)\right)=
\sum_{|h|<2K}
\left(1-\frac{|h|}{2K}\right)
\sum_{p_1,p_2\in\mathbb{P}}
e\!\left(h(f(p_1)-f(p_2))\right)\nonumber
\\
&=
\sum_{|h|<2K}
\left(1-\frac{|h|}{2K}\right)|S(h)|^2=
|S(0)|^2+
2\sum_{h=1}^{2K-1}
\left(1-\frac{h}{2K}\right)|S(h)|^2\nonumber
\\
&\ge P^2+2\sum_{h=1}^{K}
\left(1-\frac{h}{2K}\right)|S(h)|^2 \ge P^2+
\sum_{k=1}^{K}|S(k)|^2.\label{eq:fejer-prime-moment}
\end{align}
The last inequality holds because all terms are nonnegative and 
\[
2\left(1-\frac{h}{2K}\right)\ge1,
\qquad 1\le h\le K.
\]
From (\ref{cauchy}) and (\ref{eq:fejer-prime-moment}), using the nonnegativity of the kernel 
$\Phi_{2K}(t)$ and extending summation over the second variable from primes to all integers in the interval, we obtain 
\begin{align}
&V_n^2(K,x,y)
\le
K\left(
\sum_{p_1\in\mathbb{P}}\sum_{p_2\in\mathbb{P}}
\Phi_{2K}\!\left(f(p_1)-f(p_2)\right)-P^2
\right)\nonumber\\
&\le
K\left(
\sum_{p\in\mathbb{P}}\sum_{x-y<m\le x}
\Phi_{2K}\!\left(f(p)-f(m)\right)-P^2
\right)=
K\sum_{|h|<2K}
\left(1-\frac{|h|}{2K}\right)
S(h)\overline{T(h)}-KP^2\nonumber\\
&=
KP(y-P)+ 2K\operatorname{Re}
\sum_{h=1}^{2K-1}
\left(1-\frac{h}{2K}\right)
S(h)\overline{T(h)}.
\label{eq:V-fejer-primes}
\end{align}

In the last equality in~(\ref{eq:V-fejer-primes}), we used $S(0)=P$ and, since 
$x,y$ are integers, $\# \mathbb{I}=y$ and $T(0)=y$. The terms with frequencies $h$ and $-h$ are combined using the identities 
$S(-h)=\overline{S(h)}$ and 
$T(-h)=\overline{T(h)}$.

To estimate the mixed sum, observe that for each $h$
\[
2\operatorname{Re}\bigl(S(h)\overline{T(h)}\bigr)
\le |S(h)|^2+|T(h)|^2,
\]
which follows from the inequality $|S(h)-T(h)|^2\ge0$.

From identity~(\ref{eq:fejer-prime-moment}) and the inequality between the double sums in~(\ref{eq:V-fejer-primes}), taking account of the nonnegativity of the weights $1-\frac{h}{2K}$, we obtain 
\begin{align*}
&P^2+2\sum_{h=1}^{2K-1}
\left(1-\frac{h}{2K}\right)|S(h)|^2\le
Py+2\operatorname{Re}\sum_{h=1}^{2K-1}
\left(1-\frac{h}{2K}\right)
S(h)\overline{T(h)}\\
&\quad\le
Py+\sum_{h=1}^{2K-1}
\left(1-\frac{h}{2K}\right)|S(h)|^2
+\sum_{h=1}^{2K-1}
\left(1-\frac{h}{2K}\right)|T(h)|^2.
\end{align*}
Subtracting from both sides 
$P^2+\sum\limits_{h=1}^{2K-1}\Big(1-\frac{h}{2K}\Big)|S(h)|^2$, we find 
\[
\sum_{h=1}^{2K-1}
\left(1-\frac{h}{2K}\right)|S(h)|^2
\le
P(y-P)+
\sum_{h=1}^{2K-1}
\left(1-\frac{h}{2K}\right)|T(h)|^2.
\]

Consequently, 
\begin{align}\label{fejer-reduction}
V_n^2(K,x,y)
&\le
2K\sum_{h=1}^{2K-1}
\left(1-\frac{h}{2K}\right)|S(h)|^2\le
2KP(y-P)+
2K\sum_{h=1}^{2K-1}
\left(1-\frac{h}{2K}\right)|T(h)|^2
\nonumber\\
&\le
2KP(y-P)+
2K\sum_{h=1}^{2K-1}|T(h)|^2.
\end{align}
Put 
\[
D=2^{n-1},\qquad H=2K-1,\qquad
M_D(H)=\sum_{h=1}^{H}|T(h)|^D.
\]
For $n\ge 3$, apply H\"older's inequality with exponents 
\[
p=\frac{D}{2},\qquad p'=\frac{D}{D-2},
\qquad \frac{1}{p}+\frac{1}{p'}=1.
\]
We obtain 
\begin{align}\label{eq:holder-moment}
\sum_{h=1}^{H}|T(h)|^2
&=\sum_{h=1}^{H}
\left(|T(h)|^2\cdot 1\right)\le
\left(\sum_{h=1}^{H}|T(h)|^D\right)^{\frac{2}{D}}
\left(\sum_{h=1}^{H}1\right)^{1-\frac{2}{D}}=H^{1-\frac{2}{D}}M_D(H)^{\frac{2}{D}}.
\end{align}
For $n=2$, we have $D=2$, and (\ref{eq:holder-moment}) holds as an equality.

Substituting (\ref{eq:holder-moment}) into (\ref{fejer-reduction}), we find 
\[
V_n^2(K,x,y)
\le 2KP(y-P)+
2K H^{1-\frac{2}{D}}M_D(H)^{\frac{2}{D}}.
\]
Consequently, by the inequality 
$\sqrt{a+b}\le\sqrt a+\sqrt b$ for $a,b\ge0$,
\[
V_n(K,x,y)
\le \sqrt{2KP(y-P)}+
\sqrt{2K}\,H^{\frac{1}{2}-\frac{1}{D}}M_D(H)^{\frac{1}{D}}.
\]
Taking account of $H<2K$ and $\frac{1}{2}-\frac{1}{D}\ge0$, we obtain 
\begin{equation}\label{eq:V-via-moment}
V_n(K,x,y)
\le \sqrt{2KP(y-P)}+
2^{1-\frac{1}{D}}K^{1-\frac{1}{D}}M_D(H)^{\frac{1}{D}}.
\end{equation}

For each $1\le h\le H$, make the substitution $m=x-y+u$. Since $x$ and $y$ are integers, we have 
\[
T(h)=\sum_{u=1}^{y}e\bigl(\varphi_h(u)\bigr),
\qquad
\varphi_h(u)=h f(x-y+u).
\]
Applying Lemma~\ref{Weyl} with $Q=y$, $j=n-1$
 and taking account of $D=2^{n-1}$, we obtain by the triangle inequality 
\begin{align}\label{eq:weyl-applied}
|T(h)|^D
&\le
(2y)^{D-n}
\sum_{|r_1|<y}\cdots\sum_{|r_{n-1}|<y}
\left|
\sum_{u\in I_{n-1}(\mathbf r)}
e\bigl(
\Delta_{n-1}(\varphi_h(u);r_1,\ldots,r_{n-1})
\bigr)
\right|,
\end{align}
where $\mathbf r=(r_1,\ldots,r_{n-1})$, and 
$I_{n-1}(\mathbf r)\subset[1,y]$ are the intervals in Lemma~\ref{Weyl}, each containing at most $y$ integer points.

After $n-1$ successive differences, we have 
\begin{equation}\label{eq:linear-difference}
\Delta_{n-1}(\varphi_h(u);r_1,\ldots,r_{n-1})
=
n!\alpha h r_1\cdots r_{n-1}u+C_{h,\mathbf r},
\end{equation}
where $C_{h,\mathbf r}$ is independent of $u$. Indeed, each differencing step reduces the degree of the polynomial by one, and its leading coefficient is successively multiplied by 
$nr_1,(n-1)r_2,\ldots,2r_{n-1}$.

On the right-hand side of~(\ref{eq:weyl-applied}), separate the tuples for which $r_1\cdots r_{n-1}=0$. In this case, 
\[
\Delta_{n-1}(\varphi_h(u);r_1,\ldots,r_{n-1})=0,
\]
so the inner sum does not exceed $y$. The number of such tuples is bounded above by 
\[
(n-1)(2y-1)^{n-2}\ll y^{n-2}.
\]
Consequently, their total contribution to the right-hand side of~(\ref{eq:weyl-applied}) is 
\[
\ll (2y)^{D-n}y^{n-2}y
\ll y^{D-1}.
\]

Now suppose that $r_1\cdots r_{n-1}\ne0$. Since $C_{h,\mathbf r}$ is independent of $u$, (\ref{eq:linear-difference}) and Lemma~\ref{linear-sum} imply 
\[
\begin{aligned}
\left|
\sum_{u\in I_{n-1}(\mathbf r)}
e\bigl(
\Delta_{n-1}(\varphi_h(u);r_1,\ldots,r_{n-1})
\bigr)
\right|&=
\left|
\sum_{u\in I_{n-1}(\mathbf r)}
e\bigl(n!\alpha h r_1\cdots r_{n-1}u\bigr)
\right|
\\
&\le \min\left\{
y,\frac{1}{2\|n!\alpha h r_1\cdots r_{n-1}\|}
\right\}.
\end{aligned}
\]
Using the identity $\|t\|=\|-t\|$, we may restrict the resulting bound to positive shifts, absorbing the factor 
$2^{n-1}$ into the constant. Thus, 
\begin{equation}\label{eq:T-after-differencing}
|T(h)|^D
\ll y^{D-1} + y^{D-n}
\sum_{1\le r_1<y}\cdots\sum_{1\le r_{n-1}<y}
\min\left\{
y,\frac{1}{2\|n!\alpha h r_1\cdots r_{n-1}\|}
\right\}.
\end{equation}
Summing~(\ref{eq:T-after-differencing}) over $1\le h\le H$, we obtain 
\begin{equation}\label{eq:moment-after-differencing}
M_D(H)
\ll
Hy^{D-1}+
y^{D-n}
\sum_{h=1}^{H}
\sum_{1\le r_1<y}\cdots\sum_{1\le r_{n-1}<y}
\min\left\{
y,\frac{1}{2\|n!\alpha h r_1\cdots r_{n-1}\|}
\right\}.
\end{equation}
Put 
\[
Z=n!Hy^{n-1},
\qquad
F(v)=\min\left\{y,\frac{1}{2\|\alpha v\|}\right\},
\]
where, when $\|\alpha v\|=0$, we set $F(v)=y$.

Let $c(v)$ denote the number of tuples 
$(h,r_1,\ldots,r_{n-1})$ satisfying 
\[
1\le h\le H,\qquad
1\le r_i<y\quad (1\le i\le n-1),
\qquad
n!hr_1\cdots r_{n-1}=v.
\]
Since each such tuple satisfies $1\le v\le Z$, grouping terms with the same value of the product gives 
\[
\sum_{h=1}^{H}
\sum_{1\le r_1<y}\cdots
\sum_{1\le r_{n-1}<y}
F(n!hr_1\cdots r_{n-1})
=
\sum_{v\le Z}c(v)F(v).
\]

Let 
\[
\tau_n(v)=
\sum_{\substack{d_1\cdots d_n=v\\d_1,\ldots,d_n\ge1}}1
\]
be the number of representations of $v$ as a product of 
$n$ ordered positive integer factors. Each tuple counted by $c(v)$ corresponds to the factorization 
\[
v=(n!h)r_1\cdots r_{n-1}.
\]
Distinct tuples correspond to distinct ordered factorizations, so 
\[
c(v)\le\tau_n(v).
\]
Taking account also of $F(v)\ge0$, we deduce from (\ref{eq:moment-after-differencing}) that 
\begin{equation}\label{eq:moment-divisors}
M_D(H)
\ll
Hy^{D-1}
+
y^{D-n}\sum_{v\le Z}\tau_n(v)F(v).
\end{equation}
By H\"older's inequality, 
\[
\sum_{v\le Z}\tau_n(v)F(v)
\le
\left(\sum_{v\le Z}\tau_n^3(v)\right)^{\frac{1}{3}}
\left(\sum_{v\le Z}F^{\frac{3}{2}}(v)\right)^{\frac{2}{3}}.
\]
Since 
\[
F^{\frac{3}{2}}(v)=
\min\left\{y^{\frac{3}{2}},\frac{1}{(2\|\alpha v\|)^{\frac{3}{2}}}\right\},
\]
we obtain 
\begin{equation}\label{eq:moment-Holder}
M_D(H)\ll Hy^{D-1}+
y^{D-n}
\left(\sum_{v\le Z}\tau_n^3(v)\right)^{\frac{1}{3}}\left(
\sum_{v\le Z}
\min\left\{y^{\frac{3}{2}},\frac{1}{(2\|\alpha v\|)^{\frac{3}{2}}}\right\}
\right)^{\frac{2}{3}}.
\end{equation}
We estimate the first factor on the right-hand side of (\ref{eq:moment-Holder}). First, we show that 
\begin{equation}\label{eq:divisor-pointwise-holder}
\tau^3_n(v)\le\tau_{n^3}(v).
\end{equation}
For a prime $p$ and an integer $a\ge1$, we have 
\[
\tau_n(p^a)
=
\binom{a+n-1}{n-1}
=
\prod_{j=1}^{a}\frac{n+j-1}{j}.
\]
For every $j\ge1$, the following inequality holds:
\[
\left(\frac{n+j-1}{j}\right)^3
\le
\frac{n^3+j-1}{j}.
\]
Since the function $f_1(t)=t^3$ is convex for $t\ge0$, and the numbers $1-\frac{1}{j}$ and $\frac{1}{j}$ for $j\ge1$
 are nonnegative and sum to $1$, Jensen's inequality gives 
\begin{align*}
\left(\frac{n+j-1}{j}\right)^3
=
\left(
\left(1-\frac{1}{j}\right)\cdot1
+\frac{1}{j}\cdot n
\right)^3\le
\left(1-\frac{1}{j}\right)\cdot1^3
+\frac{1}{j}\cdot n^3
=
\frac{n^3+j-1}{j}.
\end{align*}
Multiplying these inequalities over $1\le j\le a$, we obtain 
\[
\tau_n^3(p^a)
\le
\prod_{j=1}^{a}\frac{n^3+j-1}{j}
=
\tau_{n^3}(p^a).
\]
By the multiplicativity of $\tau_n$ and 
$\tau_{n^3}$, this implies (\ref{eq:divisor-pointwise-holder}) for every positive integer $v$.

Next, for any integer $s\ge2$, separating the last factor in the factorization 
$v=d_1\cdots d_s$, we find

\begin{align}
\sum_{v\le Z}\tau_s(v)
&=
\sum_{\substack{d_1\cdots d_s\le Z\\d_1,\ldots,d_s\ge1}}1=
\sum_{\substack{d_1\cdots d_{s-1}\le Z\\
                d_1,\ldots,d_{s-1}\ge1}}
\left[\frac{Z}{d_1\cdots d_{s-1}}\right]\le
Z
\sum_{\substack{d_1\cdots d_{s-1}\le Z\\
                d_1,\ldots,d_{s-1}\ge1}}
\frac{1}{d_1\cdots d_{s-1}}\nonumber\\
&\le
Z\left(\sum_{d\le Z}\frac{1}{d}\right)^{s-1}
\le
Z(1+\log Z)^{s-1}.\label{eq:divisor-sum}
\end{align}

In the penultimate step, we extended the summation to all tuples 
$1\le d_1,\ldots,d_{s-1}\le Z$.

Applying (\ref{eq:divisor-sum}) with $s=n^3$ and taking account of (\ref{eq:divisor-pointwise-holder}), we obtain 
\begin{equation}\label{eq:divisor-second-moment-holder}
\sum_{v\le Z}\tau_n^3(v)
\le
\sum_{v\le Z}\tau_{n^3}(v)
\ll
Z(\log 2Z)^{n^3-1}.
\end{equation}

Substituting (\ref{eq:divisor-second-moment-holder}) into (\ref{eq:moment-Holder}), we deduce

\begin{align}\label{eq:moment-before-spacing-holder}
M_D(H)\ll Hy^{D-1}+
y^{D-n}Z^{\frac{1}{3}}\Big(\log(2Z)\Big)^{\frac{n^3-1}{3}}
\left(\sum_{v\le Z}\min\left\{y^{\frac{3}{2}},\frac{1}{(2\|\alpha v\|)^{\frac{3}{2}}}\right\}
\right)^{\frac{2}{3}}.
\end{align}
Write 
\[
\alpha=\frac{a}{q}+\frac{\theta}{q^2},
\qquad
(a,q)=1,\qquad |\theta|\le1.
\]
Put
\[
U(Z,y)=
\sum_{v\le Z}
\min\left\{y^{\frac{3}{2}},\frac{1}{(2\|\alpha v\|)^{\frac{3}{2}}}\right\}.
\]
We prove that 
\begin{equation}\label{eq:spacing-square-holder}
U(Z,y)
\ll
\left(\frac{Z}{q}+1\right)(y^{\frac{3}{2}}+q\sqrt{y}).
\end{equation}

For $q=1$, this bound follows from 
$U(Z,y)\le Zy^{\frac{3}{2}}$.

Now suppose that $q\ge2$. Partition the integers 
$1\le v\le Z$ into consecutive blocks, each containing at most 
\[
L=\left[\frac{q}{2}\right]
\]
 numbers, where $[t]$ denotes the greatest integer not exceeding $t$. Since $L\ge \frac{q}{3}$, the number of blocks does not exceed 
\[
\frac{Z}{L}+1\le\frac{3Z}{q}+1.
\]

If distinct integers $v$ and $w$ belong to the same block, then 
\[
0<|v-w|<\frac{q}{2}.
\]
By the condition $(a,q)=1$, the number $a(v-w)$
 is not divisible by $q$, so 
\[
\left\|\frac{a(v-w)}{q}\right\|\ge\frac{1}{q}.
\]
Consequently,

\begin{align*}
\|\alpha(v-w)\|
\ge
\left\|\frac{a(v-w)}{q}\right\|
-\frac{|\theta|\,|v-w|}{q^2}\ge
\frac{1}{q}-\frac{|v-w|}{q^2}
\ge\frac{1}{2q}.
\end{align*}

Thus, within each block, the points 
$\alpha v$ modulo one are pairwise separated by a distance of at least $\frac{1}{2q}$.

Consider a nonempty block $\mathcal B$. For each $v\in\mathcal B$
, put 
\[
\xi_v=\alpha v-\left[\alpha v+\frac{1}{2}\right].
\]
Then 
\[
-\frac{1}{2}\le \xi_v<\frac{1}{2},
\qquad
|\xi_v|=\|\alpha v\|.
\]
Since $\xi_v-\xi_w$ differs from $\alpha(v-w)$ by an integer, for distinct $v,w\in\mathcal B$ we have 
\[
|\xi_v-\xi_w|
\ge \|\xi_v-\xi_w\|
=\|\alpha(v-w)\|
\ge \frac{1}{2q}.
\]

Partition the numbers $\xi_v$, $v\in\mathcal B$ into nonnegative and negative ones. First, consider the nonnegative numbers. If this group is nonempty, arrange its elements in increasing order:
\[
0\le u_1<u_2<\cdots<u_s<\frac{1}{2}.
\]
By what was proved above, 
\[
u_{j+1}-u_j\ge\frac{1}{2q}
\qquad (1\le j<s).
\]
Consequently, for $2\le j\le s$
\[
u_j
=u_1+\sum_{i=1}^{j-1}(u_{i+1}-u_i)
\ge \frac{j-1}{2q},
\]
whence 
\[
\frac{1}{(2u_j)^{\frac{3}{2}}}\le\frac{q^{\frac{3}{2}}}{(j-1)^{\frac{3}{2}}}.
\]
Bounding the first term by $y^{\frac{3}{2}}$, we obtain 
\begin{align*}
\sum_{\substack{v\in\mathcal B\\ \xi_v\ge0}}
\min\left\{y^{\frac{3}{2}},\frac{1}{(2\|\alpha v\|)^{\frac{3}{2}}}\right\}=
\sum_{j=1}^{s}
\min\left\{y^{\frac{3}{2}},\frac{1}{(2u_j)^{\frac{3}{2}}}\right\}&\le
y^{\frac{3}{2}}+\sum_{j=2}^{s}
\min\left\{y^{\frac{3}{2}},\frac{q^{\frac{3}{2}}}{(j-1)^{\frac{3}{2}}}\right\}\\
&\le
y^{\frac{3}{2}}+\sum_{r=1}^{\infty}
\min\left\{y^{\frac{3}{2}},\frac{q^{\frac{3}{2}}}{r^{\frac{3}{2}}}\right\}.
\end{align*}
If there are no nonnegative numbers $\xi_v$, the left-hand side is zero and the bound remains valid.

For the negative numbers $\xi_v$, consider the positive numbers $-\xi_v$. The distance between any two of them is still at least $\frac{1}{2q}$, since 
\[
|(-\xi_v)-(-\xi_w)|
=|\xi_v-\xi_w|
\ge\frac{1}{2q}.
\]
Applying the same argument to these numbers, we find 
\[
\sum_{\substack{v\in\mathcal B\\ \xi_v<0}}
\min\left\{y^{\frac{3}{2}},\frac{1}{(2\|\alpha v\|)^{\frac{3}{2}}}\right\}
\le
y^{\frac{3}{2}}+\sum_{r=1}^{\infty}
\min\left\{y^{\frac{3}{2}},\frac{q^{\frac{3}{2}}}{r^{\frac{3}{2}}}\right\}.
\]
Thus, 
\[
\sum_{v\in\mathcal B}
\min\left\{y^{\frac{3}{2}},\frac{1}{(2\|\alpha v\|)^{\frac{3}{2}}}\right\}
\le
2y^{\frac{3}{2}}+2\sum_{r=1}^{\infty}
\min\left\{y^{\frac{3}{2}},\frac{q^{\frac{3}{2}}}{r^{\frac{3}{2}}}\right\}.
\]

The function 
\[
g(u)=\min\left\{y^{\frac{3}{2}},\frac{q^{\frac{3}{2}}}{u^{\frac{3}{2}}}\right\},
\qquad u>0,
\]
is nonnegative and nonincreasing. Hence 
\[
g(r)\le\int\limits_{r-1}^{r}g(u)\,du
\qquad (r\ge1),
\]
and consequently, 
\[
\sum_{r=1}^{\infty}
\min\left\{y^{\frac{3}{2}},\frac{q^{\frac{3}{2}}}{r^{\frac{3}{2}}}\right\}
\le
\int\limits_0^\infty
\min\left\{y^{\frac{3}{2}},\frac{q^{\frac{3}{2}}}{u^{\frac{3}{2}}}\right\}du=\int\limits_0^{\frac{q}{y}}y^{\frac{3}{2}}du+
\int\limits_{\frac{q}{y}}^{\infty}\frac{q^{\frac{3}{2}}}{u^{\frac{3}{2}}}du=q\sqrt{y}+2q\sqrt{y}=3q\sqrt{y}.
\]
As a result, for each block $\mathcal B$ we obtain 
\[
\sum_{v\in\mathcal B}
\min\left\{y^{\frac{3}{2}},\frac{1}{(2\|\alpha v\|)^{\frac{3}{2}}}\right\}
\le 2y^{{\frac{3}{2}}}+6q\sqrt{y}.
\]

Taking account of the fact that the number of blocks does not exceed $\frac{3Z}{q}+1$, we finally find 
\[
U(Z,y)
\le
\left(\frac{3Z}{q}+1\right)(2y^{{\frac{3}{2}}}+6q\sqrt{y})
\ll
\left(\frac{Z}{q}+1\right)(y^{{\frac{3}{2}}}+q\sqrt{y}),
\]
which proves the required bound (\ref{eq:spacing-square-holder}).

From (\ref{eq:moment-before-spacing-holder}) and (\ref{eq:spacing-square-holder}), it follows that 
\begin{equation*}
M_D(H)\ll 
Hy^{D-1}+y^{D-n}Z^{\frac{1}{3}}\Big(\log(2Z)\Big)^{\frac{n^3-1}{3}}\left[
\left(\frac{Z}{q}+1\right)(y^{{\frac{3}{2}}}+q\sqrt{y})
\right]^{\frac{2}{3}}.
\end{equation*}

Observe that 
\[
Z\left(\frac{Z}{q}+1\right)
\left(y^{\frac{3}{2}}+q\sqrt{y}\right)
=
Z^2y^{\frac{3}{2}}
\left(
\frac{1}{q}+\frac{1}{y}+\frac{1}{Z}+\frac{q}{Zy}
\right).
\]

Therefore, 
\begin{equation*}
  M_D(H)\ll Hy^{D-1}+ y^{D-n}Zy\Bigl(\log(2Z)\Bigr)^{\frac{n^3-1}{3}}
\left(\frac{1}{q}+\frac{1}{y}+\frac{1}{Z}+\frac{q}{Zy}
\right)^{\frac{2}{3}}.
\end{equation*}

Since $Z=n!Hy^{n-1}$ and $H=2K-1<2y$, we have 
\begin{equation}\label{eq:Z-relations}
\begin{gathered}
y^{D-n}Zy=n!Hy^D, \qquad \frac{1}{Z}\le\frac{1}{y},
\qquad \frac{q}{Zy}\le\frac{q}{Hy^n}, \\
Z<2n!y^n, \qquad \log(2Z)\ll\log(2y).
\end{gathered}
\end{equation}

Introducing the notation 
\[
\Delta_H=\frac{1}{q}+\frac{1}{y}+\frac{q}{Hy^n},
\]
we find 
\[
M_D(H)
\ll
Hy^{D-1}
+
Hy^D\Delta_H^{\frac{2}{3}}
\Big(\log(2y)\Big)^{\frac{n^3-1}{3}}.
\]
Since $\Delta_H\ge \frac{1}{y}$ and $y\ge2$, we have 
\[
Hy^{D-1}
\le
Hy^D\Delta_H^{\frac{2}{3}}
\Big(\log(2y)\Big)^{\frac{n^3-1}{3}}.
\]
Thus, 
\begin{equation}\label{eq:moment-final-holder}
M_D(H) \ll Hy^D\left(
\frac{1}{q}+\frac{1}{y}+\frac{q}{Hy^n}
\right)^{\frac{2}{3}}
\Big(\log(2y)\Big)^{\frac{n^3-1}{3}}.
\end{equation}

Substituting (\ref{eq:moment-final-holder}) into (\ref{eq:V-via-moment}), we obtain 
\begin{equation}\label{MDHHolder}
  V_n(K,x,y)\ll \sqrt{2KP(y-P)}+K^{1-\frac{1}{D}}H^{\frac{1}{D}}y
\left(\frac{1}{q}+\frac{1}{y}+\frac{q}{Hy^n}
\right)^{\frac{2}{3D}}
\Big(\log(2y)\Big)^{\frac{n^3-1}{3D}}.
\end{equation}

We now apply the Cauchy--Schwarz inequality to the sum in (\ref{eq:moment-divisors}):
\[
\left(\sum_{v\le Z}\tau_n(v)F(v)\right)^2
\le
\left(\sum_{v\le Z}\tau_n^2(v)\right)
\left(\sum_{v\le Z}F^2(v)\right).
\]
Since 
\[
F^2(v)=
\min\left\{y^2,\frac{1}{4\|\alpha v\|^2}\right\},
\]
we obtain 
\begin{equation}\label{eq:moment-cauchy}
M_D(H)\ll Hy^{D-1}+ y^{D-n}
\left(\sum_{v\le Z}\tau_n^2(v)\right)^{1/2}\left(
\sum_{v\le Z}
\min\left\{y^2,\frac{1}{4\|\alpha v\|^2}\right\}
\right)^{1/2}.
\end{equation}
We estimate the first factor on the right-hand side of (\ref{eq:moment-cauchy}). First, we show that 
\begin{equation}\label{eq:divisor-pointwise}
\tau^2_n(v)\le\tau_{n^2}(v).
\end{equation}
For a prime $p$ and an integer $a\ge1$, we have 
\[
\tau_n(p^a)
=
\binom{a+n-1}{n-1}
=
\prod_{j=1}^{a}\frac{n+j-1}{j}.
\]
For every $j\ge1$, the following inequality holds:
\[
\left(\frac{n+j-1}{j}\right)^2 \le
\frac{n^2+j-1}{j},
\]
since 
\[
j(n^2+j-1)-(n+j-1)^2 =
(j-1)(n-1)^2\ge 0.
\]
Multiplying these inequalities over $1\le j\le a$, we obtain 
\[
\tau_n^2(p^a)\le
\prod_{j=1}^{a}\frac{n^2+j-1}{j}
= \tau_{n^2}(p^a).
\]
By the multiplicativity of $\tau_n$ and 
$\tau_{n^2}$, this implies (\ref{eq:divisor-pointwise}) for every positive integer $v$.

Applying (\ref{eq:divisor-sum}) with $s=n^2$ and taking account of (\ref{eq:divisor-pointwise}), we obtain 
\begin{equation}\label{eq:divisor-second-moment}
\sum_{v\le Z}\tau_n^2(v)
\le
\sum_{v\le Z}\tau_{n^2}(v)
\ll
Z(\log 2Z)^{n^2-1}.
\end{equation}

Substituting (\ref{eq:divisor-second-moment}) into (\ref{eq:moment-cauchy}), we deduce

\begin{align}\label{eq:moment-before-spacing}
M_D(H)\ll &
Hy^{D-1}+y^{D-n}Z^{\frac{1}{2}}\Big(\log(2Z)\Big)^{\frac{n^2-1}{2}}
\left(
\sum_{v\le Z}
\min\left\{y^2,\frac{1}{4\|\alpha v\|^2}\right\}
\right)^{\frac{1}{2}}.
\end{align}
Since $0\le F(v)\le y$, we have 
\[
F^2(v)=F^{\frac{3}{2}}(v)\sqrt{F(v)}
\le \sqrt{y}\,F^{\frac{3}{2}}(v).
\]
Therefore, (\ref{eq:spacing-square-holder}) gives 
\begin{equation}\label{eq:spacing-square}
\sum_{v\le Z}F^2(v)
\le
\sqrt{y}\sum_{v\le Z}F^{\frac{3}{2}}(v)
\ll
\left(\frac{Z}{q}+1\right)(y^2+qy).
\end{equation}

From (\ref{eq:moment-before-spacing}) and (\ref{eq:spacing-square}), it follows that 
\begin{equation*}
M_D(H)\ll Hy^{D-1}+y^{D-n}Z^{\frac{1}{2}}\Big(\log(2Z)\Big)^{\frac{n^2-1}{2}}\left[
\left(\frac{Z}{q}+1\right)(y^2+qy)
\right]^{\frac{1}{2}}.
\end{equation*}

Observe that 
\[
Z\left(\frac{Z}{q}+1\right)(y^2+qy)=
Z^2y^2 \left(\frac{1}{q}+\frac{1}{y}+\frac{1}{Z}+\frac{q}{Zy}
\right).
\]
Therefore, 
\begin{equation*}
  M_D(H)\ll Hy^{D-1}+ y^{D-n}Zy\Bigl(\log(2Z)\Bigr)^{\frac{n^2-1}{2}}
\left(\frac{1}{q}+\frac{1}{y}+\frac{1}{Z}+\frac{q}{Zy}
\right)^{\frac{1}{2}}.
\end{equation*}

Using relations (\ref{eq:Z-relations}), we obtain 
\begin{align*}
M_D(H)
&\ll
Hy^{D-1}+
Hy^D\Delta_H^{\frac{1}{2}}
\bigl(\log(2y)\bigr)^{\frac{n^2-1}{2}}.
\end{align*}
Since $\Delta_H\ge\frac{1}{y}$ and 
$\log(2y)>1$, we have 
\[
Hy^{D-1}\le
Hy^D\Delta_H^{\frac{1}{2}}
\bigl(\log(2y)\bigr)^{\frac{n^2-1}{2}}.
\]
Consequently, 
\begin{align}\label{eq:moment-cauchy-simplified}
M_D(H)\ll Hy^D\Delta_H^{\frac{1}{2}}
\bigl(\log(2y)\bigr)^{\frac{n^2-1}{2}}.
\end{align}
Substituting the bound (\ref{eq:moment-cauchy-simplified}) into (\ref{eq:V-via-moment}), we find 
\begin{align}\label{eq:V-cauchy-moment}
V_n(K,x,y)\ll
\sqrt{2KP(y-P)}+K^{1-\frac{1}{D}}H^{\frac{1}{D}}y
\Delta_H^{\frac{1}{2D}}
\bigl(\log(2y)\bigr)^{\frac{n^2-1}{2D}}.
\end{align}

Finally, since $K\le H<2K$, we have 
\[
K^{1-\frac{1}{D}}H^{\frac{1}{D}}
\le 2^{\frac{1}{D}}K, \qquad
\Delta_H\le\Delta.
\]
From (\ref{MDHHolder}) and (\ref{eq:V-cauchy-moment}), taking account of $D=2^{n-1}$, we obtain 
\begin{align}\label{eq:V-final1}
V_n(K,x,y)\ll 
\sqrt{2KP(y-P)}+Ky\min\left\{
\Delta^{\frac{1}{2^n}}
\bigl(\log(2y)\bigr)^{\frac{n^2-1}{2^n}},\Delta^{\frac{1}{3\cdot 2^{n-2}}}
\bigl(\log(2y)\bigr)^{\frac{n^3-1}{3\cdot 2^{n-1}}}
\right\}.
\end{align}
By the Brun--Titchmarsh inequality~\cite[p.~121, Theorem~2]{MontgomeryVaughan1973}, we have 
\[
P=\pi(x)-\pi(x-y)\ll\frac{y}{\log(2y)}
\]
and taking account of $0\le P\le y$, we have 
\[
\sqrt{2KP(y-P)}
\le \sqrt{2KPy}
\ll \frac{y\sqrt K}{\sqrt{\log(2y)}}
=\frac{Ky}{\sqrt{K\log(2y)}}.
\]
Consequently, we finally obtain 
\begin{align}\label{eq:V-final}
V_n(K,x,y)
\ll Ky\left[
\frac{1}{\sqrt{K\log(2y)}}+\min\left\{
\Delta^{\frac{1}{2^n}}\bigl(\log(2y)\bigr)^{\frac{n^2-1}{2^n}},\Delta^{\frac{1}{3\cdot 2^{n-2}}}
\bigl(\log(2y)\bigr)^{\frac{n^3-1}{3\cdot 2^{n-1}}}\right\}\right].
\end{align}
This proves Theorem~\ref{theo1}.

\section{Proof of Corollary~\ref{cor1}}
The conditions on $y$ and $q$ imply 
\[
\Delta=
\frac{1}{q}+\frac{1}{y}+\frac{q}{Ky^n}
\ll (\log(2y))^{-B}.
\]
Moreover, the condition on $K$ gives 
\[
\frac{1}{\sqrt{K\log(2y)}}
\ll (\log(2y))^{-A}.
\]

If $B=2^nA+n^2-1$, then 
\[
\Delta^{\frac{1}{2^n}}
\bigl(\log(2y)\bigr)^{\frac{n^2-1}{2^n}}
\ll (\log(2y))^{-A}.
\]
Otherwise, $B=3\cdot 2^{n-2}A+\frac{n^3-1}{2},$
and hence 
\[
\Delta^{\frac{1}{3\cdot 2^{n-2}}}
\bigl(\log(2y)\bigr)^{\frac{n^3-1}{3\cdot 2^{n-1}}}
\ll (\log(2y))^{-A}.
\]
Thus, at least one of the expressions under the minimum sign in (\ref{eq:V-final}) satisfies the required bound. Applying (\ref{eq:V-final}) proves the corollary.

\section{Proof of Corollary~ \ref{cor2}}
All arguments in the proof of the theorem preceding the application of the Brun--Titchmarsh inequality remain valid when summation over primes is replaced by summation over all integers in $\mathbb I$. In this case, 
\[
S(h)=T(h),\qquad P=\#\mathbb I=y.
\]
Consequently, 
\[
\sqrt{2KP(y-P)}=0.
\]
The bound (\ref{eq:V-final1}), obtained before applying the Brun--Titchmarsh inequality, proves the corollary.

\end{document}